\documentclass[oneside]{article}
\usepackage[utf8]{inputenc}
\usepackage[english]{babel}  

\usepackage{amsthm}
\newtheorem{theorem}{Theorem}[section]
\newtheorem{lemma}{Lemma}[section]
\newtheorem{proposition}{Proposition}[section]

\theoremstyle{definition}

\usepackage{amsmath}
\usepackage{amsfonts}
\usepackage{amssymb}
\usepackage{graphicx}
\numberwithin{equation}{section}
\usepackage{lipsum}
\usepackage{mathtools}
\usepackage{scalefnt}
\usepackage{empheq}
\usepackage[a4paper,left=3.5cm,right=3cm,top=3cm,bottom=3cm]{geometry}
\usepackage{makeidx}
\usepackage{bm}

\makeindex

\usepackage{indentfirst}

\begin{document}
	
	\title{Existence for noncoercive nonlinear elliptic equations with two lower-order terms}
	
	\author{Thi Tam Dang\thanks{Department of Mathematics, Leopold Franzens Universit\"at Innsbruck, Austria, e-mail: Thi-Tam.Dang@uibk.ac.at}, Trung Hau Hoang \thanks{Department of Mathematics, Leopold Franzens Universit\"at Innsbruck, Austria, e-mail: Trung-Hau.Hoang@uibk.ac.at} \\}
	
	\maketitle
	
	\begin{abstract}
	This paper considers a class of noncoercive nonlinear elliptic problems with coefficients defined in Marcinkiewicz and Lorentz spaces. We prove the existence of a solution for the corresponding Dirichlet problem and investigate the higher integrability properties of the solution. 
	\end{abstract}

	\section{Introduction}
This paper considers a general noncoercive nonlinear elliptic problem of the form:
\begin{equation}\label{1.1}
	\begin{cases}
		&-\mathrm{div} A(x,u,\nabla u) + B(x,u,\nabla u) + G(x,u) = F \quad \text{in} \quad \Omega,\\
		& u = 0 \quad \text{ on} \quad \partial \Omega,
	\end{cases}
\end{equation}
where $\Omega$ is a bounded open subset of $\mathbb{R}^{N}$ with $N \ge 2$ and $F \in W^{-1, p^{\prime}}(\Omega)$. The operator $A(x, z, \xi): \Omega \times \mathbb{R} \times \mathbb{R}^{N} \rightarrow \mathbb{R}^{N} $ is a Carath\'eodory vector field which meets the following conditions: for a.e. $x\in \Omega$, $z \in \mathbb{R}$ and $\xi, \xi_{*} \in \mathbb{R}^{N}$
\begin{align}
	&\left\langle A(x,z, \xi), \xi \right\rangle \ge \alpha \vert \xi \vert^{p}- \big( b(x) \vert z \vert \big)^{p}- \varphi(x)^{p}, \quad \alpha >0, \label{1.2}\\
	& \vert A(x,z, \xi ) \vert \le \beta \vert \xi \vert^{p-1} + \big( b(x) \vert z \vert \big)^{p-1} +  \varphi(x)^{p-1}, \quad \beta >0, \label{1.3} \\
	& \left\langle A(x,z, \xi) - A(x,z, \xi_{*}), \xi - \xi_{*} \right\rangle  >0  \qquad \text{for} \ \ \xi \ne \xi_{*}, \label{1.4}
\end{align}
where $b(x) \in L^{N, \infty}(\Omega)$, and $\varphi(x) \in L^{p}(\Omega)$ with $ 1<p<N$.\\
The operator $B(x,z, \xi) : \Omega \times \mathbb{R} \times \mathbb{R}^{N} \rightarrow \mathbb{R}$ is a Carath\'eodory function satisfies
\begin{align}\label{1.5}
	& \vert B(x,z, \xi) \vert \le c(x) \vert \xi \vert^{p-1} + \phi(x),
\end{align}
for a.e. $x\in \Omega$, $z \in \mathbb{R}$, $\xi, \xi_{*} \in \mathbb{R}^{N}$, the coefficient $c(x) \in L^{N,1}(\Omega)$ and $\phi \in L^1(\Omega)$.\\
Finally, the vector field $G(x,z): \Omega \times \mathbb{R} \rightarrow \mathbb{R}$ is a Carath\'eodory function satisfying
\begin{align}
	& G(x,z)z \ge 0, \label{1.6}\\
	& \vert G(x,z)\vert  \le d(x)\vert z \vert^{\lambda}+ \psi(x), 	\ \ 0 \le \lambda < \frac{N(p-1)}{N-p}, \label{1.7}
\end{align}
where $d(x) \in L^{s^{\prime},1}(\Omega)$ with $ s = \frac{N(p-1)}{N-p} \frac{1}{\lambda}, \frac{1}{s} + \frac{1}{s^{\prime}} =1$ and $\psi(x) \in L^1(\Omega)$.

A simple example of the model problem \eqref{1.1} can be read as follows:
\begin{equation}
	\begin{cases}
		&-\triangle_{p}u + \mathrm{div}(K(x) \vert u \vert^{p-2} u) +c(x) \vert \nabla u \vert^{p-1}= -\mathrm{div} \left( \frac{x}{\vert x \vert^{N- \nu}}\right)  \quad \text{in} \ \ \Omega,      \\
		& u=0 \quad \text{ on} \ \ \partial \Omega,
	\end{cases}
\end{equation}
for any $\frac{N}{2}- 1 < \nu \le N-1$. For a.e. $x \in \Omega$, we define the vector field $K: \Omega \rightarrow \mathbb{R}^{N}$ is a mesurable function satisfying $\vert K(x) \vert \le b(x)^{p-1}$ with $b(x) \in L^{N, \infty}(\Omega)$.
Thus, the operator $A$ is a combination of the following two terms:
\begin{align}
	A(x,z, \xi):= -\triangle_{p}u + \mathrm{div}(K(x) \vert u \vert^{p-2} u).
\end{align}

In the case the operator $B$ and $G$ vanishes, the problem \eqref{1.1} has been studied in \cite{Farroni } where the authors proved the existence of the solution for the quasilinear elliptic problems with singluarity in the lower order term. The
existence of the renormalized solutions has been studied in \cite{Guibe} in which the operator $A$ is independent of $u$ with no singluarity behavior. 

The problem \eqref{1.1} has two main features: first, the vector field $A$ introduces a singular lower-order term influenced by $u$, described by the coefficient $b(x)$. The property that the coefficient $b(x)$ belongs to the Marcinkiewicz spaces $L^{N, \infty}(\Omega)$ leads to a lack of compactness of the operator $K(x,u)$. That is, the operator $K(x,u)$ does not belong to $L^{p^{\prime}}(\Omega)$ in general, and the term $\mathrm{div} (K(x,u))$ is not an element of the dual space $W^{-1, p^{\prime}}(\Omega)$. Second, the operator $B$ is characterized by a coefficient, denoted $c(x)$, belonging to the Lorentz spaces $L^{N,1}(\Omega)$, which exhibit unboundedness. Typically, the presence of these singularity and unboundedness properties causes the problem to be noncoercive. 


For any $\varepsilon \in [0,1]$, we set
\begin{align}\label{1.11}
	\zeta_{\varepsilon}(x)= \varepsilon c(x) + (1-\varepsilon)b(x), \quad \forall x \in \Omega.
\end{align}
We further define $X_{\varepsilon}(\Omega)$ to be a convex subset of $L^{N, q}(\Omega)$, with $1 \le q \le \infty$ such that for $	\zeta_{\varepsilon}(x) \in X_{\varepsilon}(\Omega)$ satisfies 
\begin{align}\label{1.12}
	\| \zeta_{\varepsilon} - T_{k}\zeta_{\varepsilon} \|_{N,q} \le \theta, \quad \theta = \frac{1}{S_{N,p}} \left(\frac{\alpha}{2p} \right)^{\frac{1}{p}}, 
\end{align}
for any $k>0$, $S_{N,p}$ denotes the Sobolev constant given by Theorem \eqref{Th: 2.1} below. \\
Our main existence result is stated in the next theorem.
\begin{theorem}\label{Th: 1.1}
	Let the assumptions \eqref{1.2}-\eqref{1.7} be fulfilled. Furthermore, we assume that 
	\begin{align}\label{1.13}
		\zeta_{\varepsilon} \in X_{\varepsilon}(\Omega).
	\end{align}
	Then the problem \eqref{1.1} admits a solution.
\end{theorem}

The novelty of this study lies in establishing the existence theorem by an approximation method. To achieve this, a sequence of approximation problems is considered (see Section~3.1), and the existence of a solution to these problems is demonstrated using the Leray--Schauder fixed point theorem. The primary challenge is deriving a priori estimates for the sequence of solutions, which is facilitated by fulfilling our assumption \eqref{1.13}. Subsequently, we establish the compactness of the sequence of solutions to the approximation problems through rigorous testing with an admissible test function. Then we show that the limit of the sequence is a solution to the original problem (see Section~3.2).

Consider for a moment the problem \eqref{1.1} in the linear case and the coefficient $c=0$, Stampacchia's results (see $\cite{Stampacchia}$) proved that $u \in L^{r^{*}}(\Omega)$ by assuming the coefficients of the lower order terms in appropriate Lebesgue spaces. For nonlinear elliptic problems, in \cite{Boca} Stampacchia's results were extended by assuming $b(x)\in L^{N/p-1}(\Omega)$. In $\cite{GreMos}$ a similar result was achieved for nonlinear elliptic problems with lower order term coefficients in Marcinkiewicz spaces. Their analysis assumed that these coefficients are in $L^{\infty}(\Omega)$. In addition, in a recent publication \cite{Farroni }, the authors extended the previous results to noncoercive quasilinear elliptic operators with a singular lower order term.

The intriguing question that arises is whether the aforementioned results can be extended to a broader class of noncoercive nonlinear elliptic problem, involving coefficients of lower order terms in both the Marcinkiewicz and Lorentz spaces. In our present work we aim to extend the above results to our problem. More precisely, we will prove the higher summability of a solution under the appropriate properties of the given data.

The structure of the paper is as follows: we start with the Preliminaries, where we recall the definition and some properties of Lorentz spaces, which play a fundamental role in our analytical framework. In Section~3 we will demonstrate the existence of a solution by an approximation procedure. We will show that the sequence of solutions to this approximation problem converges to a solution of the original problem. Finally, Section~4 investigates the regularity properties of a solution.


\section{Preliminaries}
This section recalls a definition and representative properties of Lorentz spaces, which are used extensively throughout this paper. For further information on Lorentz spaces, we refer the readers to $\cite{Lorentz}$.

Let $\Omega \subset \mathbb{R}^{N}$ be a bounded domain. The distribution function of $f$ is given by
\begin{align}\label{2.1}
	\omega_{f}(h) =  \vert E_{h} \vert = \vert \left\lbrace x \in \Omega: \vert f(x) \vert > h \right\rbrace \vert,
\end{align}
where $ \vert E_{h}\vert$ is the Lebesgue measure of $E_{h}$. For given $1 <r < \infty$ and $1 \le s < \infty$,  the Lorentz spaces denoted by $L^{r,s}(\Omega)$ consist of all measurable functions $f$ defined on $\Omega$ which satisfy 
\begin{align}
	\| f \|_{r,s}^{r} =  \int_0^{\infty} (\omega_{f}(h))^{\frac{s}{r}} h^{s-1} dh < \infty.
\end{align}
We remark that  $L^{r,s} (\Omega)$ becomes Banach space when it is endowed with the norm $\| \cdot \|_{r,s}$. In the case $r=s$, the Lorentz space simplifies to the Lebesgue space $L^{r}(\Omega)$. If $s= \infty$, the class $L^{r, \infty}(\Omega)$ consists of all measurable functions $f$ defined on $\Omega$ which satisfy 
\begin{align}
	\| f \|_{r, \infty}^{r} = \sup_{t>0} h^{r} \omega_{f}(h) <  \infty,
\end{align}
and it coincides with the Marcinkiewicz class, the weak-$L^{r}(\Omega)$. For Lorentz spaces, the following continuous inclusions apply:
\begin{align}
	L^{q}(\Omega) \subset L^{r,s}(\Omega) \subset L^{r,q}(\Omega) \subset L^{r,\infty}(\Omega) \subset L^{s}(\Omega)
\end{align}
for $ 1 \le s <r< q \le  \infty$. Furthermore, for any $f \in L^{r,s}(\Omega)$, $g \in L^{r^{\prime}, s^{\prime}}(\Omega)$, where $1 <r < \infty$, $1 \le s \le \infty$, $\frac{1}{r}+ \frac{1}{r^{\prime}} = 1$, $\frac{1}{s} + \frac{1}{s^{\prime}} =1$, the following generalized H\"older inequality holds
\begin{align}\label{2.4}
	\| f g \|_1 \le  \| f\|_{r,s} \| g \|_{ r^{\prime}, s^{\prime}}.
\end{align}
The distance of a function $f \in L^{r, \infty}(\Omega)$ to $L^{\infty}(\Omega)$ can be chracterized as follows: 
 \begin{align}
 	\lim_{k \rightarrow \infty} \| f -T_{k} f \|_{r, \infty} = 	\mathrm{dist}_{ L^{r, \infty}(\Omega)  } (f, L^{\infty}(\Omega)),
 \end{align}
where $T_{k}(s): \mathbb{R} \rightarrow \mathbb{R}$ denotes the truncation operator at level $\pm k$, i.e. $T_{k}(s) = \frac{s}{\vert s \vert} \min \left\lbrace \vert s \vert, k \right\rbrace $. For more about the distance to $L^{\infty}(\Omega)$ and some applications, see $\cite{Carozza}$.

 The Sobolev embedding theorem in the framework of Lorentz spaces is stated in the next theorem.
\begin{theorem}\label{Th: 2.1}(\cite{Betta, Farroni , Lorentz})
		Let us assume that $1<p<N, 1\le q \le p$. Then, any function $g\in W_0^{1,1}(\Omega)$ satisfying $\vert\nabla g \vert \in L^{p,q}(\Omega)$ belongs to $L^{p^{\star},p}(\Omega)$ where $p^{*}= \frac{Np}{N-p}$ and
	\begin{align}\label{2.6}
		\Vert g \Vert_{p^{*},q} \le S_{N,p} \Vert \nabla g\Vert_{p,q},
	\end{align}
	where $S_{N,p}= \omega_{N}^{-1/N} \frac{p}{N-p}$ and $\omega_{N}$ stands for the measure of the unit ball in $\mathbb{R}^{N}$.
\end{theorem}
We need the following the weak compactness lemma for the proof of the Theorem \ref{Th: 1.1}.
\begin{lemma}\label{Le: 2.1} (Lemma~3.1, \cite{Farroni })
	Let $\mathcal{B}$ be a nonempty subset of $W_0^{1,p}(\Omega)$. Assume that there exists a constant $C>0$ such that
	\begin{align}
		\| \nabla u \|_{L^{p}( \Omega \backslash E_{k})}^{p} \le C (1 + \| u \|_{L^{p}(\Omega \backslash E_{k})}^{p}),
	\end{align}
	for any $k>0$ and $u \in \mathcal{B}$, where $E_{k} := \left\lbrace x \in \Omega: \vert u(x) \vert \ge k\right\rbrace $. Then, there exists a constant $M>0$ such that
	\begin{align}
		\| u \|_{W^{1,p}(\Omega)} \le M,
	\end{align}
	for any $u \in \mathcal{B}$.
\end{lemma}
We conclude this section with the well known Leray--Schauder fixed point theorem.
\begin{theorem}(Leray--Schauder fixed point) (\cite{Farroni , Gilbarg}) \label{Th: 2.2}
	Let $\mathcal{F}$ be a compact mapping of a Banach space $X$ into itself, and all suppose there exists a constant $M$ such that $\| x \| <M$ for all $x \in X$ and $t \in [0,1]$ satisfying $x = t \mathcal{F}(x)$. Then, $\mathcal{F}$ has a fixed point.	
\end{theorem}

	
\section{Proof of theorem}	
We first consider the approximation problems, then showing that the approximation problems have a solution via the Leray-Schauder fixed-point theorem that converges to a solution of our original problem \eqref{1.1}.

\subsection{The approximation problems}		
For each $n \in \mathbb{N}$, for any $\varepsilon \in [0,1]$, we set
\begin{align*}
	\vartheta_{n}^{\varepsilon}(x) := \frac{T_{n}\zeta_{\varepsilon}(x)}{\zeta_{\varepsilon}(x)}, \quad \text{a.e.} \ \ x \in \Omega,
\end{align*}
where $\zeta_{\varepsilon}$ is given by \eqref{1.11}. We consider the following approximating problems:
\begin{align}\label{3.1}
	- \mathrm{div}A_{n}(x,u_{n}, \nabla u_{n}) + B_{n}(x,u_{n},\nabla u_{n}) + G_{n}(x,u_{n}) = F,
\end{align}
where the vector field $A_{n}: (x,z, \xi) \in \Omega \times \mathbb{R} \times \mathbb{R}^{N} \rightarrow \mathbb{R}^{N}$ is defined by
\begin{align*}
	A_{n}(x,z,\xi) = A(x, \vartheta_{n}^{\varepsilon}z, \xi).
\end{align*}
The vector field $A_{n}$ shares properties with $A$ as follows:
\begin{align}
&	\left\langle A_{n}(x,z,\xi), \xi \right\rangle \ge \alpha \vert \xi \vert^{p} - \big( T_{n} \zeta_{0}(x) \vert z \vert \big)^{p} - \varphi(x)^{p}, \label{3.2}\\
& \vert A_{n}(x,z,\xi) \vert \le \beta \vert \xi \vert^{p-1} +  \big( T_{n} \zeta_{0}(x) \vert z \vert \big)^{p-1} + \varphi(x)^{p-1}, \label{3.3}\\
& \left\langle  A_{n}(x,z,\xi) -  A_{n}(x,z, \xi_{*}), \xi - \xi_{*}    \right\rangle >0  \qquad \text{for} \ \ \xi \ne \xi_{*}. \label{3.4}
\end{align}
Moreover, we let $	B_{n}(x,z,\xi):  \Omega \times \mathbb{R} \times \mathbb{R}^{N} \rightarrow \mathbb{R}^{N} $ satisfy
\begin{align}\label{3.5}
	\vert B_{n}(x,u,\xi) \vert \le T_{n}\zeta_{1}(x)\vert \xi \vert^{p-1} + \phi(x).
\end{align}	
Finally, the vector field $G_{n}(x,z) = T_{n}G(x,u)$ fulfils the following conditions:
\begin{align}
	& G_{n}(x,z) z \ge 0,\label{3.6} \\
	&\vert G_{n}(x,z) \vert \le d(x) \vert z\vert^{\lambda} + \psi(x). \label{3.7}
\end{align}	
Our goal is to find $u_{n} \in W_0^{1,p}(\Omega)$, which solves the approximation problem \eqref{3.1}. The existence of the solution $u_{n}$ can be proven by applying the Leray-Schauder fixed point theorem. 
For this, we need the following technical lemmas, which provide an a priori estimate for $u_{n}$. The next lemma gives an a priori estimate for $ \| \vert \nabla S_{k} u_{n} \vert^{p-1} \|_{N^{\prime}, \infty} $ where $S_{k}(u_{n})$ is defined by \eqref{3.11} below.
\begin{lemma} \label{Le: 3.1}
	Let $\Omega$ be an open subset of $\mathbb{R}^{N}$ with finite measure and that $1<p<N$. Let $u$ be a measurable function satisfying $T_{k}(u_{n}) \in W_0^{1,p}(\Omega)$, for every positive $k>0$, and such that
	\begin{align}
		\int_{\Omega} \vert \nabla T_{\sigma}(S_{k}u_{n})\vert^{p} dx \le Mk + L, \quad \forall k>0,
	\end{align}
	where $M$ and $L$ are given constants. Then $ \vert \nabla S_{k}u_{n} \vert^{p-1}$ belongs to $L^{N^{\prime}, \infty}(\Omega)$ and
	\begin{align}\label{3.9}
		\| \vert \nabla S_{k}u_{n}\vert^{p-1} \|_{N^{\prime}, \infty} \le 2 C(N,p) \left[ \frac{ \| \phi\|_1}{C_1}  + \vert \Omega \vert^{\frac{1}{N^{\prime}} - \frac{1}{p^{\prime}} } L^{ \frac{1}{p^{\prime}}  }    \right]
	\end{align}	
	where $C(N,p)$ is a constant depending only on $N$ and $p$.	
\end{lemma}

\begin{proof}
	For $k>0$, for all $s \in \mathbb{R}$, the remainder $S_{k}(s)$ of the truncation $T_{k}(s)$ is given by
	\begin{align}
		S_{k}(s) = s -T_{k}(s),
	\end{align}
	that is 
	\begin{align} S_{k}(s)= \label{3.11}
		\begin{cases}
			0, \quad \vert s \vert \le k, \\
			(\vert s \vert - k)\mathrm{sign}(s), \quad  \vert s \vert >k.
		\end{cases}
	\end{align}	
	For fixed $\sigma >0$, using $T_{\sigma}(S_{k}(u_{n}))$ as a test function of the approximation equation \eqref{3.1}, we obtain
	\begin{equation}\label{3.12}
		\begin{aligned}
			& \int_{\Omega} A_{n}(x, u_{n}, \nabla u_{n}) \cdot \nabla T_{\sigma}(S_{k}(u_{n})) dx + \int_{\Omega} B_{n}(x,u_{n}, \nabla u_{n}) T_{\sigma}(S_{k}(u_{n})) dx \\
			& \qquad + \int_{\Omega} G_{n}(x,u_{n}) T_{\sigma}(S_{k}(u_{n})) dx = \int_{\Omega} F T_{\sigma}(S_{k}(u_{n})) dx. 
		\end{aligned}
	\end{equation}
	By the definition of $S_{k}(u_{n})$ given by \eqref{3.11} and the ellipticity condition \eqref{3.2} on $A_{n}$, we get
	\begin{equation}\label{3.13}
		\begin{aligned}
			& \int_{\Omega} A_{n}(x, u_{n}, \nabla u_{n}) \cdot \nabla T_{\sigma}(S_{k}(u_{n})) dx \\
			& \quad = \int_{ \left\lbrace  k \le\vert u_{n} \vert \le k+\sigma \right\rbrace } A_{n}(x, u_{n}, \nabla u_{n}) \cdot \nabla u_{n} dx \\
			& \quad \ge \int_{ \left\lbrace  k \le\vert u_{n} \vert \le k+\sigma \right\rbrace } \Big(  \alpha \vert \nabla u_{n} \vert^{p} - ( T_{n} (\zeta_0) \vert u_{n} \vert)^{p} - \varphi(x)^{p} \Big)  dx\\
			& \quad \ge \alpha \int_{\Omega} \vert \nabla T_{\sigma}(S_{k}(u_{n})) \vert^{p} dx- \int_{ \left\lbrace  k \le\vert u_{n} \vert \le k+\sigma \right\rbrace } (T_{n}(\zeta_0) \vert u_{n} \vert)^{p}  dx - \| \varphi \|_{p}^{p}.
		\end{aligned}
	\end{equation}
	For any $n \ge m$, where $m$ is a positive integer to be chosen later, we have
	\begin{align*}
		T_{n} \zeta_{\varepsilon} \le T_{m} \zeta_{\varepsilon} + (\zeta_{\varepsilon} - T_{m} \zeta_{\varepsilon}), \quad \varepsilon = 0,1,
	\end{align*}
	which implies 
	\begin{equation}
		\begin{aligned}
			& \int_{ \left\lbrace  k \le\vert u_{n} \vert \le k+\sigma \right\rbrace } (T_{n}(\zeta_0)  u_{n} )^{p}  dx = \| T_{n} \zeta_0 \vert u_{n} \vert \chi_{ \left\lbrace k \le \vert u_{n} \vert \le k + \sigma \right\rbrace }\|_{p}^{p}\\
			&\quad  \le  \| T_{m}( \zeta_0) u_{n} \chi_{ \left\lbrace k \le \vert u_{n} \vert \le k + \sigma \right\rbrace } \|_{p}^{p} + \|  (\zeta_0-T_{m} (\zeta_0))  u_{n} \chi_{ \left\lbrace k \le \vert u_{n} \vert \le k + \sigma \right\rbrace } \|_{p}^{p}.
		\end{aligned}	
	\end{equation}
	Using the H\"older inequality \eqref{2.4}, the generalized Sobolev embedding Theorem \ref{Th: 2.1}, and the assumption \eqref{1.13}, we deduce that
	\begin{equation}
		\begin{aligned}
			\| T_{m} (\zeta_0 ) u_{n}  \chi_{ \left\lbrace k \le \vert u_{n} \vert \le k + \sigma \right\rbrace } \|_{p}^{p} &\le m^{p}
			\|  u_{n}  \chi_{ \left\lbrace k \le \vert u_{n} \vert \le  k + \sigma \right\rbrace } \|_{p}^{p}\\
			& \le m^{p}  \| 1 \|_{N,\infty}^{p} \|  u_{n} \chi_{ \left\lbrace k \le \vert u_{n} \vert \le  k + \sigma \right\rbrace }  \|_{ p^{*},p }^{p}  \\
			& \le m^{p}  \vert \Omega \vert^{\frac{p}{N}}S_{N,p}^{p}  \|  \nabla T_{\sigma}(S_{k}(u_{n})) \|_{p}^{p},
		\end{aligned}
	\end{equation}
	and
	\begin{equation}
		\begin{aligned}
			\|  (\zeta_0-T_{m} (\zeta_0) ) u_{n} \chi_{ \left\lbrace k \le \vert u_{n} \vert \le k + \sigma \right\rbrace } \|_{p}^{p} &\le \| \zeta_0-T_{m} (\zeta_0) \|_{N, \infty}^{p}  \| u_{n} \chi_{ \left\lbrace k \le \vert u_{n} \vert \le  k + \sigma \right\rbrace }    \|_{ p^{*},p }^{p}  \\
			& \le S_{N,p}^{p} \| \zeta_0-T_{m} \zeta_0 \|_{N, \infty}^{p}  \|  \nabla T_{\sigma}(S_{k}(u_{n})) \|_{p}^{p}\\
			& \le \frac{\alpha}{2p}  \|  \nabla T_{\sigma}(S_{k}(u_{n})) \|_{p}^{p}.
		\end{aligned}
	\end{equation}
	By the definition of $S_{k}(u_{n}) $, $S_{k}(u_{n})=0$ for $\vert u_{n}\vert \le k$, and the growth assumption \eqref{3.5} on $B_{n}$, we get
	\begin{equation}\label{3.17}
		\begin{aligned}
			&	\left| \int_{\Omega} B_{n}(x, u_{n}, \nabla u_{n}) T_{\sigma}(S_{k}(u_{n}) ) dx \right| \\
			&  \le \sigma \int_{\Omega} \vert 	B_{n}(x, u_{n}, \nabla u_{n}) \vert dx \\
			& \le \sigma \left[ \int_{ \left\lbrace \vert u_{n} \vert > k\right\rbrace  } T_{n}(\zeta_1 ) \vert \nabla u_{n} \vert^{p-1} dx + \int_{\Omega} \phi(x) dx \right] \\
			& \le  \sigma \int_{ \left\lbrace \vert u_{n} \vert > k\right\rbrace  } T_{n} (\zeta_1) \vert \nabla u_{n} \vert^{p-1} dx + \sigma \| \phi\|_1. 
		\end{aligned}
	\end{equation}
	For the first term on the right side of \eqref{3.17}, we use the H\"older inequality and the assumption \eqref{1.13} to obtain
	\begin{equation}\label{3.18}
		\begin{aligned}
			&\sigma \int_{ \left\lbrace \vert u_{n} \vert > k\right\rbrace  } T_{n}( \zeta_1) \vert \nabla u_{n} \vert^{p-1} dx \\
			& \quad \le \sigma  \int_{ \left\lbrace \vert u_{n} \vert > k\right\rbrace  } T_{m} (\zeta_1 )\vert \nabla u_{n} \vert^{p-1} dx + \sigma  \int_{ \left\lbrace \vert u_{n} \vert > k\right\rbrace  } (\zeta_1 - T_{m} (\zeta_1))  \vert \nabla u_{n} \vert^{p-1} dx\\
			& \quad \le \sigma m \int_{\Omega} \vert \nabla u_{n} \chi_{ \left\lbrace \vert u_{n} \vert > k\right\rbrace } \vert^{p-1} dx + \sigma \| \zeta_1 - T_{m}(\zeta_1) \|_{N,1} \| \vert \nabla S_{k}(u_{n}) \vert^{p-1} \|_{N^{\prime}, \infty} \\
			&\quad  \le \sigma m  \| 1 \|_{N,1} \| \vert \nabla S_{k}(u_{n})\|_{N^{\prime}, \infty} + \sigma \theta \| \vert \nabla S_{k}(u_{n}) \vert^{p-1} \|_{N^{\prime}, \infty} \\
			&\quad  \le\sigma ( m N \vert \Omega \vert^{\frac{1}{N}} +  \theta) \| \vert \nabla S_{k}(u_{n}) \vert^{p-1} \|_{N^{\prime}, \infty},
		\end{aligned}
	\end{equation}
	where $\theta$ is given by \eqref{1.13}. The sign condition \eqref{3.6} on $G_{n}$ leads to
	\begin{equation}\label{3.19}
		\begin{aligned}
			\int_{\Omega} G_{n}(x,u_{n}) T_{\sigma}(S_{k}(u_{n})) dx \ge 0.
		\end{aligned}
	\end{equation}
	Applying Young's inequality yields 
	\begin{equation}\label{3.10}
		\begin{aligned}
			\int_{\Omega} F  T_{\sigma}(S_{k}(u_{n})) dx &\le \| F \| \| \nabla T_{\sigma}(S_{k}(u_{n})) \|_{p} \\
			& \le \frac{\alpha}{2p} \| \nabla T_{\sigma}(S_{k}(u_{n})) \|_{p}^{p} + \frac{2^{p^{\prime}/p}}{p^{\prime}\alpha^{p^{\prime}/p} } \| F \|^{p^{\prime}}_{p^{\prime}}.
		\end{aligned}
	\end{equation}
	We set
	\begin{align}
		C_1 := \frac{\alpha}{p^{\prime}} - m^{p}\vert \Omega \vert^{\frac{p}{N}} S_{N,p}^{p}.
	\end{align}
	In view of \eqref{3.13}-\eqref{3.10}, we obtain
	\begin{equation}
		\begin{aligned}
			\|  \nabla T_{\sigma}(S_{k}(u_{n})) \|_{p}^{p} \le M \sigma + L, \qquad \forall \sigma >0,
		\end{aligned}
	\end{equation}
	where 
	\begin{equation}\label{eq: 3.53}
		\begin{aligned}
			&M = \frac{1}{C_1} \left[  \Big(m N \vert \Omega \vert^{\frac{1}{N}} + \theta \Big) \| \vert \nabla S_{k}(u_{n}) \vert^{p-1} \|_{N^{\prime},\infty } + \| \phi\|_1 \right] , \\
			& L =  \frac{1}{C_1}  \Big( \| \varphi \|_{p}^{p} + \frac{2^{p^{\prime}/p}}{p^{\prime}\alpha^{p^{\prime}/p} } \| F \|^{p^{\prime}}_{p^{\prime}} \Big).	
		\end{aligned}
	\end{equation}
	By Lemma~4.1 of \cite{Guibe}, one has
	\begin{equation}
		\begin{aligned}
			& \| \vert  \nabla S_{k}(u_{n}) \vert^{p-1} \|_{N^{\prime},\infty } \\
			& \quad \le C(N,p) \left[ M + \vert \Omega \vert^{\frac{1}{N^{\prime}} - \frac{1}{p^{\prime}} } L^{ \frac{1}{p^{\prime}}  }    \right] \\
			&\quad \le C(N,p) \left[  \frac{mN \vert \Omega \vert^{\frac{1}{N}} + \theta}{C_1}   \| \vert \nabla S_{k}(u_{n}) \vert^{p-1} \|_{N^{\prime},\infty }  + \frac{ \| \phi\|_1}{C_1}  + \vert \Omega \vert^{\frac{1}{N^{\prime}} - \frac{1}{p^{\prime}} } L^{ \frac{1}{p^{\prime}}  }    \right]. 
		\end{aligned}
	\end{equation}
	We can choose $m$ to be large enough to guarantee that
	\begin{align}
		C(N,p) \frac{mN \vert \Omega \vert^{\frac{1}{N}} + \theta}{C_1}  \le \frac{1}{2}.
	\end{align}
	Therefore we have
	\begin{equation}\label{eq: 3.56}
		\begin{aligned}
			\| \vert \nabla S_{k}(u_{n}) \vert^{p-1} \|_{N^{\prime},\infty } \le 2 C(N,p) \left[ \frac{ \| \phi\|_1}{C_1}  + \vert \Omega \vert^{\frac{1}{N^{\prime}} - \frac{1}{p^{\prime}} } L^{ \frac{1}{p^{\prime}}  }    \right],
		\end{aligned}
	\end{equation}
	which completes the proof.		
\end{proof}
\begin{lemma}\label{Le: 3.2}
	Let the assumptions of the Theorem \ref{Th: 1.1} be fulfilled.  Let $u_{n}$ be a measurable function satisfying $T_{k}(u_{n}) \in W_0^{1,p}(\Omega)$, for every positive $k>0$, there exist a constant $C$ such that
	\begin{align}
		\| \nabla u \|_{L^{p}( \Omega \backslash E_{k})}^{p} \le C (1 + \| u \|_{L^{p}(\Omega \backslash E_{k})}^{p}),
	\end{align}
	where $E_{k}$ are defined by \eqref{2.1}.
\end{lemma}
\begin{proof}
For $k>0$, testing the equation \eqref{3.1} with $T_{k}(u_{n})$, we get
	\begin{equation}\label{3.27}
		\begin{aligned}
			&	\int_{\Omega} A_{n}(x,u_{n}, \nabla u_{n}) \cdot \nabla T_{k}(u_{n} )dx + \int_{\Omega} B_{n}(x,u_{n},\nabla u_{n}) T_{k}(u_{n})dx\\
			& \qquad + \int_{\Omega} G_{n}(x,u_{n}) T_{k}(u_{n})dx = \int_{\Omega} F T_{k}(u_{n})dx.
		\end{aligned}
	\end{equation}
	The ellipticity condition \eqref{3.2} on $A_{n}$ implies
	\begin{equation}\label{3.28}
		\begin{aligned}
			&\int_{\Omega} A_{n}(x,u_{n},\nabla u_{n}) \cdot \nabla T_{k}(u_{n}) dx\\
			& \quad = \int_{ \left\lbrace \vert u_{n} \vert \le k\right\rbrace   }A_{n}(x,u_{n},\nabla u_{n}) \cdot \nabla u_{n} dx \\
			& \quad \ge \int_{ \left\lbrace \vert u_{n} \vert \le k\right\rbrace } \Big(  \alpha \vert \nabla u_{n} \vert^{p} - \big( T_{n} \zeta_{0} \vert u_{n} \vert \big)^{p} - \varphi(x)^{p} \Big) dx\\
			& \quad \ge \alpha \int_{\Omega} \vert \nabla T_{k}(u_{n}) \vert^{p}dx - \int_{ \left\lbrace \vert u_{n} \vert \le k\right\rbrace } \Big( T_{n} (\zeta_0 )\vert u_{n} \vert \Big)^{p} dx - \| \varphi \|_{p}^{p}.
		\end{aligned}
	\end{equation}
	With the help of the H\"older inequality \eqref{2.4}, the Soblolev embedding theorem \eqref{Th: 2.1}, and the assumption \eqref{1.13}, we gain
	\begin{equation}
		\begin{aligned}
			& \int_{ \left\lbrace \vert u_{n} \vert \le k\right\rbrace } \Big( T_{n} (\zeta_0 )\vert u_{n} \vert \Big)^{p} dx = \| (T_{n}(\zeta_{0}) ) u_{n} \chi_{ \left\lbrace \vert u_{n} \vert \le k\right\rbrace } \|_{p}^{p} \\
			&	\quad \le  	\| (T_{m} (\zeta_{0})) u_{n} \chi_{ \left\lbrace \vert u_{n} \vert \le k\right\rbrace } \|_{p}^{p} + \| ( \zeta_0 - T_{m} (\zeta_0)) u_{n} \chi_{ \left\lbrace \vert u_{n} \vert \le k\right\rbrace } \|_{p}^{p}\\
			& \quad \le m^{p} \| u_{n} \|_{ L^{p}(\Omega \backslash E_{k} )}^{p} +  \| \zeta_0 - T_{m} (\zeta_0) \|_{N, \infty}^{p} \| T_{k} u_{n} \|_{p^{*},p}^{p} \\
			& \quad \le 	m^{p} \| u_{n} \|_{ L^{p}(\Omega \backslash E_{k} )}^{p}+ \| \zeta_0 - T_{m}(\zeta_0)  \|_{N, \infty}^{p} S_{N,p}^{p} \| \nabla T_{k} u_{n} \|_{p}^{p}\\
			& \quad \le 	m^{p} \| u_{n} \|_{ L^{p}(\Omega \backslash E_{k} )}^{p} + \frac{\alpha}{2p} \| \nabla T_{k} u_{n} \|_{p}^{p}.
		\end{aligned}
	\end{equation}
	By the definition of $S_{k}(u_{n})$ \eqref{3.11}, the growth assumption \eqref{3.5} on $B_{n}$, and the H\"older inequality, one obtains
	\begin{equation}\label{3.30}
		\begin{aligned}
			&\left| \int_{\Omega} B_{n}(x, u_{n}, \nabla u_{n}) T_{k}(u_{n})dx  \right| \\
			& \quad \le k \int_{\Omega} \vert  B_{n}(x, u_{n}, \nabla u_{n}) \vert dx \\
			& \quad \le k \left[ \int_{ \left\lbrace \vert u_{n} \vert \le k\right\rbrace } T_{n}( \zeta_1 ) \cdot \vert \nabla u_{n} \vert^{p-1} dx + \int_{ \left\lbrace \vert u_{n} \vert > k\right\rbrace } T_{n}( \zeta_1 ) \cdot \vert \nabla u_{n} \vert^{p-1} dx + \int_{\Omega} \phi(x) \right] \\
			&\quad  \le \frac{2^{p/p^{\prime}} k^{p}}{p \alpha^{p/p^{\prime}}} \| T_{n}(\zeta_1  )\|_{p}^{p} + \frac{\alpha}{2p^{\prime}} \| \nabla T_{k}(u_{n}) \|_{p}^{p} + k \| \phi\|_{1}\\
			& \qquad + k \left[ \int_{ \left\lbrace \vert u_{n} \vert > k\right\rbrace } T_{m}(\zeta_1 ) \cdot \vert \nabla u_{n} \vert^{p-1} dx + \int_{ \left\lbrace \vert u_{n} \vert > k\right\rbrace }  \big( \zeta_1 - T_{m}( \zeta_1) \big) \cdot   \vert \nabla u_{n} \vert^{p-1} dx \right] \\
			&\quad  \le \frac{2^{p/p^{\prime}} k^{p}}{p \alpha^{p/p^{\prime}}} \| T_{n}(\zeta_1 \|_{p}^{p} + \frac{\alpha}{2p^{\prime}} \| \nabla T_{k}(u_{n}) \|_{p}^{p} + k \| \phi\|_{1}\\
			& \qquad + km N \vert \Omega\vert^{\frac{1}{N}} \| \vert \nabla \mathcal{S}_{k}u_{n}\vert^{p-1} \|_{N^{\prime}, \infty} + k \theta \vert \nabla \mathcal{S}_{k}u_{n}\vert^{p-1} \|_{N^{\prime}, \infty} \\
			& \quad \le \frac{2^{p/p^{\prime}} k^{p}}{p \alpha^{p/p^{\prime}}} \| T_{n}(\zeta_1) \|_{p}^{p} + \frac{\alpha}{2p^{\prime}} \| \nabla T_{k}(u_{n}) \|_{p}^{p} \\
			& \qquad + k \left[  (m N \vert \Omega\vert^{\frac{1}{N}}+ \theta) \| \vert \nabla \mathcal{S}_{k}u_{n}\vert^{p-1} \|_{N^{\prime}, \infty} + \| \phi \|_1 \right] ,
		\end{aligned}
	\end{equation}
	where $\theta$ is given by \eqref{1.13}. The sign condition \eqref{3.5} on $G_{n}$ leads to
	\begin{align}\label{3.31}
		\int_{\Omega} G_{n}(x, u_{n}) T_{k}(u_{n}) dx \ge 0.
	\end{align}
	According to Young's inequality, we have the following
	\begin{equation}\label{3.32}
		\begin{aligned}
			\int_{\Omega} F T_{k}(u_{n})dx &\le \| F \| \| \nabla T_{k}(u_{n}) \|_{p} \\
			& \le   \frac{2^{p^{\prime}/p}}{p^{\prime}\alpha^{p^{\prime}/p} } \| F \|^{p^{\prime}}_{p^{\prime}}+ \frac{\alpha}{2p} \| \nabla T_{k}(u_{n}) \|_{p}^{p}.
		\end{aligned}
	\end{equation}
	Therefore, using \eqref{3.28}-\eqref{3.32}, we obtain
	\begin{equation}\label{3.33}
		\begin{aligned}
			\frac{\alpha}{2 p^{\prime}}	\int_{\Omega} \vert \nabla T_{k}(u_{n})\vert^{p} dx &\le m^{p} \| u_{n} \|_{ L^{p}(\Omega \backslash E_{k} )}^{p} + \frac{2^{p/p^{\prime}} k^{p}}{p \alpha^{p/p^{\prime}}} \| T_{n}(\zeta_1) \|_{p}^{p}   +  \frac{2^{p^{\prime}/p}}{p^{\prime}\alpha^{p^{\prime}/p} } \| F \|^{p^{\prime}}_{p^{\prime}}\\
			& \quad + k \left[ (m N \vert \Omega\vert^{\frac{1}{N}}+ \theta) \| \nabla S_{k}(u_{n}) \vert^{p-1}\|_{N^{\prime}, \infty} + \| \phi \|_1    \right],
		\end{aligned}
	\end{equation}
	where the estimate of $ \| \nabla S_{k}(u_{n})\vert^{p-1}\|_{N^{\prime}, \infty }$ is given by \eqref{3.9}. This completes the proof.
\end{proof}
The existence of a solution $u_{n}\in W_0^{1,p}(\Omega)$, which solves \eqref{3.1} can be directly achieved from the following Proposition.
\begin{proposition}\label{Pro: 3.1}
	Let \eqref{1.2}-\eqref{1.7} and \eqref{1.13} be in charge. The the approximation problem \eqref{3.1} admits a solution $u_{n} \in W_0^{1,p}$.
\end{proposition}
\begin{proof}
	For any fixed $F \in W^{-1, p^{\prime}}(\Omega)$, let $\mathcal{F}: v_{n} \in W_0^{1,p}(\Omega) \mapsto u_{n} \in W_0^{1,p}(\Omega)$ be a mapping that takes $v_{n}$ to the unique solution of the problem \eqref{3.1} is compact (see more Corollary~4.1, \cite{Farroni }). It is obvious that a fixed point of $\mathcal{F}$ is a solution of the approximation \eqref{3.1}. To apply the Leray-Schauder fixed point, we need to obtain an a priori estimate for the solution $u_{n}$ by applying directly Lemma \ref{Le: 2.1}. A priori estimate for $u_{n}$ can be obtained by first establishing a priori estimate for $\vert \nabla S_{k}(u_{n}) \vert^{p-1}$ (see Lemma \ref{Le: 3.1}), where the bound depends on both $k$ and the given data. We will then derive a priori estimate for $ \nabla T_{k}(u_{n})$, also depending on $k$ and the data (see Lemma \ref{Le: 3.2}). By Lemma \ref{Le: 2.1}, we obtain
	\begin{align}\label{3.34}
		\| u_{n}\| \le M,
		\end{align}
	where $M$ is independent of $n$. The proof is concluded by applying the Leray--Schauder fixed point theorem (see Theorem \ref{Th: 2.2}).
\end{proof}



We are now ready to proceed passing to the limit to conclude the proof of the Theorem \ref{Th: 1.1}.
\subsection{Passing to the limit}
The estimate of $u_{n}$ given in \eqref{3.34} allows us to obtain
\begin{equation}\label{eq: 3.89}
	\begin{aligned}
	&	u_{n} \rightharpoonup u \qquad \text{in} \ \ W_0^{1,p}(\Omega) \ \ \text{weakly},\\
	& u_{n} \rightarrow u \qquad \text{in} \ \ L^{q}(\Omega) \ \ \text{strongly for any } \ \ q < p^{*},\\
	\end{aligned}
\end{equation}
for some $u \in W_0^{1,p}(\Omega)$. We make use of $\gamma(u_{n}-u)$ (for simplicity, $\gamma(t) := \arctan t$) as a test function for approximation problem \eqref{3.1}, we get
\begin{equation}
	\begin{aligned}
	&	\int_{\Omega} A_{n}(x, u_{n}, \nabla u_{n}) \cdot \nabla \gamma(u_{n}-u) dx + \int_{\Omega} B_{n}(x, u_{n}, \nabla u_{n}) \gamma (u_{n}-u) dx \\
	& \qquad + \int_{\Omega} G_{n}(x, u_{n}) \gamma (u_{n}-u) dx = \int_{\Omega} F \gamma(u_{n}-u) dx,
	\end{aligned}
\end{equation}
where $\nabla \gamma(u_{n} -u) = \gamma^{\prime}(u_{n}-u) (\nabla u_{n} - \nabla u)$. Since $\gamma(0)=0$, we have
\begin{align*}
	\gamma(u_{n}-u) \rightharpoonup 0 \qquad \text{in} \ \ W_0^{1,p}(\Omega) \ \ \text{weakly},
\end{align*}
which implies that
\begin{align*}
	\lim_{n \rightarrow \infty} \int_{\Omega} A_{n}(x, u_{n}, \nabla u_{n}) \cdot \nabla \gamma(u_{n}-u) dx =0.
\end{align*}
Arguing as in  the proof of \cite{Farroni }, we obtain
\begin{equation}
	\begin{aligned}
		\nabla u_{n} \rightarrow \nabla u \qquad \text{a.e. in } \ \ \Omega
 	\end{aligned}
\end{equation}
and
\begin{equation}\label{3.38}
	\begin{aligned}
		A_{n}(x, u_{n}, \nabla u_{n}) \rightharpoonup A(x,u, \nabla u) \quad \text{in} \ \ L^{p^{\prime}}(\Omega, \mathbb{R}^{N}) \ \ \text{weakly}.
	\end{aligned}
\end{equation}
%
By the estimate \eqref{3.9}, Lemma \ref{Le: 3.1} of $\vert \nabla S_{k}(u_{n})\vert^{p-1}$ together with the estimate \eqref{3.33} of $\nabla T_{k}(u_{n})$, we arrive at 
\begin{align*}
	\| \vert \nabla u_{n} \vert^{p-1} \|_{N^{\prime}, \infty} &\le \| \vert \nabla T_{k}(u_{n}) \vert^{p-1} \|_{N^{\prime}, \infty}+ \| \vert \nabla S_{k}(u_{n})\vert^{p-1} \|_{N^{\prime}, \infty}\\
	& \le C \| \vert \nabla T_{k}(u_{n}) \vert^{p-1} \|_{p} + \| \vert \nabla S_{k}(u_{n})\vert^{p-1} \|_{N^{\prime}, \infty} \le C.
\end{align*}
This results in
\begin{equation}\label{3.39}
	\begin{aligned}
	\| B_{n} (x, u_{n}, \nabla u_{n}) \|_{1} &= \int_{\Omega} \vert B_{n} (x, u_{n}, \nabla u_{n} \vert dx \\
& 	\le \int_{\Omega} T_{n}(\zeta_1) \vert \nabla u_{n} \vert^{p-1} dx + \int_{\Omega} \phi(x)dx \\
	&\le \| T_{n}\zeta_1 \|_{N,1} \| \vert \nabla u_{n} \vert^{p-1} \|_{N^{\prime}, \infty} + \| \phi \|_1 \le C.
	\end{aligned}
\end{equation}
From \eqref{3.39} and the definition of $B_{n}$, we deduce that
\begin{equation}
	\begin{aligned}
		B_{n}(x, u_{n}, \nabla u_{n}) \rightarrow B(x,u, \nabla u) \quad \text{a.e. in} \ \ \Omega.
	\end{aligned}
\end{equation}
 For every measurable set $E \subset \Omega$, we have
\begin{equation}
	\begin{aligned}
		\int_{E} \vert B_{n} (x, u_{n}, \nabla u_{n}) \vert dx &\le \| T_{n}(\zeta_1)\|_{ L^{N,1}(E) } \| \vert \nabla u_{n} \|^{p-1} \|_{L^{N^{\prime}, \infty}(E) } + \| \phi \|_{L^1(E)}\\
		& \le \| T_{n}(\zeta_0)\|_{ L^{N,1}(E) } \|  +  \| \phi \|_{L^1(E)} \le C,
	\end{aligned}
\end{equation}
which implies that
\begin{align}
	B_{n}(x, u_{n}, \nabla u_{n}) \text{ is equi-integrable.}
\end{align}
By employing the Vitali theorem, we can conclude that
\begin{align}\label{3.43}
	B_{n}(x, u_{n}, \nabla u_{n}) \rightarrow B(x, u, \nabla u) \quad \text{in } \ \ L^1(\Omega) \ \ \text{strongly}.
\end{align}
By the growth condition \eqref{3.5} on $G_{n}$ and generalized H\"older inequality, we get
\begin{equation}
	\begin{aligned}
		\| G_{n}(x, u_{n}) \|_1 = \int_{\Omega} \vert G_{n}(x, u_{n}) \vert dx & \le \int_{\Omega} d(x) \vert u_{n} \vert^{\lambda} + \psi(x) dx\\
		& \le \| d \|_{\lambda^{\prime},1} \| \vert u_{n} \vert^{\lambda} \|_{\lambda, \infty}+ \| \psi \|_1 \le C.
	\end{aligned}
\end{equation}
Similarly, we get
\begin{align}\label{3.45}
	G_{n}(x, u_{n}) \rightarrow G(x,u) \quad \text{in } \ \ L^1(\Omega) \ \ \text{strongly}.
\end{align}
In view of \eqref{3.38}, \eqref{3.43}, and \eqref{3.45}, we conclude that $u$ is the solution of \eqref{1.1}.


\section{Regularity of the solution}
In this section, we study the regularity properties of the solution of the problem \eqref{1.1}.  We generalize the result of \cite{Farroni , GreMos}. We show that $u\in L^{r^{*}}(\Omega)$ follows under appropriate assumptions on the given data.
\begin{theorem}
	Let us assume that the assumptions \eqref{1.2}-\eqref{1.7} hold with $\varphi \in L^{r}(\omega)$, $\phi \in L^{r}(\omega)$ and $F \in W^{-1, \frac{r}{p-1}}(\Omega)$ for $1<p<r<N$. For any $\zeta_{\varepsilon}$ defined in \eqref{1.11}, there exists $\eta = \eta (\alpha, N, p,r)>0$ such that when
	\begin{align}
		\mathrm{dist}_{L^{N, q}} (\zeta_{\varepsilon}, L^{q}(\Omega)) \le \eta, \quad 1 \le q \le \infty, 
	\end{align}
then any solution $u \in W_0^{1,p}(\Omega)$ satisfies
\begin{align}
	\vert u\vert^{r^{*}/p^{*}} \in W_0^{1,p}(\Omega).
\end{align}	
More precisely $u \in L^{r^{*}}(\Omega)$.
\end{theorem}

\begin{proof}
For $ F \in  W^{-1, \frac{r}{p-1}}(\Omega)$, we can take
\begin{align*}
	F = \mathrm{div} (\vert H \vert^{p-1}H).
\end{align*}
For fixed $k>0$, using in \eqref{1.1} the test function $ S_{k}(u)$ defined in \eqref{3.11}, we get
\begin{equation}
	\begin{aligned}
		\alpha \int_{E_{k}} \vert \nabla u \vert^{p} dx \le \int_{E_{k}} \Big( \zeta_0^{p} \vert u \vert^{p} + \varphi^{p} \Big) dx + \int_{E_{k}} \Big( \zeta_1 \vert \nabla u \vert^{p-1} + \phi \Big)dx + \int_{E_{k}} \vert H \vert^{p-1} \vert \nabla u \vert dx,	
	\end{aligned}
\end{equation}	
where $E_{k}$ are defined by \eqref{2.1}, $\zeta_{\varepsilon}, \varepsilon \in [0,1]$ is given by \eqref{1.11}. 

 The application of the Young's inequality gives 
\begin{equation}\label{4.4}
	\begin{aligned}
		\alpha  \int_{E_{k}} \vert \nabla u \vert^{p} dx &\le \int_{E_{k}} \Big( \zeta_0^{p} \vert u \vert^{p} + \varphi^{p} \Big) dx + \int_{E_{k}} \phi dx \\
		& \qquad+ \frac{2^{p^{\prime}/p}   }{\alpha^{p^{\prime}/p  } p^{\prime}   } \int_{E_{k}} \vert H \vert^{p} dx  + \frac{ \alpha}{2p} \int_{E_{k}} \vert \nabla u \vert^{p} dx \\
		& \qquad + \frac{ 2^{p/p^{\prime} }  }{ \alpha^{p/p^{\prime} } p } \int_{E_{k}} \zeta_1^{p}dx + \frac{\alpha}{2 p^{\prime}} \int_{E_{k}} \vert \nabla u \vert^{p}dx.\\
	\end{aligned}
\end{equation} 
Rearranging the terms on the right-hand side of \eqref{4.4}, we get
\begin{equation}\label{4.5}
	\begin{aligned}
		\frac{\alpha}{2} \int_{E_{k}}  \vert \nabla u \vert^{p} dx &\le \int_{E_{k}} \left[  M_{*}  \Big( \vert H \vert^{p}+  \zeta_1^{p} \Big) + \zeta_0^{p} \vert u \vert^{p}+ \varphi^{p} + \phi^{p} \right] dx
	\end{aligned}
\end{equation}
with $M_{*} = \max \left\lbrace  \frac{2^{p^{\prime}/p}   }{\alpha^{p^{\prime}/p  } p^{\prime}   },  \frac{ 2^{p/p^{\prime} }  }{ \alpha^{p/p^{\prime} } p }   \right\rbrace $. 

We first multiply both sides of \eqref{4.5} by $k^{p \delta-1}$. Then, for any fixed $K>0$, we integrate over the interval $[0,K]$ with respect to $k$ to get
\begin{equation}
	\begin{aligned}
		\frac{\alpha}{2} \int_{\Omega} \vert \nabla u \vert^{p} \vert T_{K}(u) \vert^{p \delta} dx \le \int_{\Omega}  \left[  M_{*}  \Big( \vert H \vert^{p}+  c^{p} \Big) + b^{p} \vert u \vert^{p}+ \varphi^{p} + \phi^{p} \right] \vert T_{K}(u) \vert^{p \delta } dx.
	\end{aligned}
\end{equation}
This can be rewritten as follows:
\begin{equation}\label{eq: 4.7}
	\begin{aligned}
		\left( \frac{\alpha}{2} \right)^{\frac{1}{p}} \| \nabla u \vert T_{K}(u) \vert^{\delta} \|_{p} &\le M_{*} \| H  \vert T_{K}(u) \vert^{\delta} \|_{p}  + M_{*} \| \zeta_1 \vert T_{K}(u) \vert^{\delta} \|_{p} + \| \zeta_0 u \vert T_{K}(u) \vert^{\delta} \|_{p} \\
		& \qquad + \| \varphi \vert T_{K}(u) \vert^{\delta} \|_{p}  +  \| \phi \vert T_{K}(u) \vert^{\delta} \|_{p}.
	\end{aligned}
\end{equation}
For $M >0$, we have
\begin{equation}
	\begin{aligned}
		\| \zeta_1 \vert T_{K}(u) \vert^{\delta} \|_{p} \le \| (\zeta_1 -T_{M}\zeta_1) \vert T_{K}(u) \vert^{\delta} \|_{p} + M \| \vert T_{K}(u) \vert^{\delta}\|_{p}.
	\end{aligned}
\end{equation}
Moreover, we have  $L^{N,1}(\Omega) \subset L^{N, \infty}(\Omega)$, thus $\zeta_1 \in L^{N, \infty}(\Omega)$. By using H\"older inequality and Sobolev embedding theorem, we get
\begin{equation}\label{eq: 4.9}
	\begin{aligned}
		\| (\zeta_1 -T_{M}(\zeta_1)) \vert T_{K}(u) \vert^{\gamma} \|_{p} & \le \| \zeta_1- T_{M}(\zeta_1) \|_{N, \infty} \| \vert T_{K}(u) \vert^{\delta} \|_{ p^{*},p}\\
		& \le  \| \zeta_1- T_{M}(\zeta_1) \|_{N, \infty} S_{N,p} \delta \| \vert \nabla u \vert \vert T_{K}(u) \vert^{\delta} \|_{p}\\
		& \le N \| \zeta_1- T_{M}(\zeta_1) \|_{N,1} S_{N,p} \delta \| \vert \nabla u \vert \vert T_{K}(u) \vert^{\delta} \|_{p}.
	\end{aligned}
\end{equation}
Here we applied the following 
\begin{align*}
	\| \zeta_{1}- T_{M}(\zeta_1)\|_{N, \infty} \le N \| \zeta_1-T_{M}(\zeta_1) \|_{N, 1}.
\end{align*}
Similarly, for $\| b u \vert T_{K}(u) \vert^{\delta} \|_{p} $,  we have
\begin{equation}\label{eq: 4.10}
	 \begin{aligned}
	 	\| \zeta_0 u \vert T_{K}(u) \vert^{\delta} \|_{p} & \le \| \zeta_1- T_{M}(\zeta_1) \|_{N, \infty} S_{N,p} (1+\delta) \| \vert \nabla u \vert \vert T_{K}(u) \vert^{\gamma} \|_{p}+ M \| u \vert T_{K} \vert^{\gamma} \|_{p}.
	 \end{aligned}
\end{equation}
We assume that
\begin{align}
	\| \zeta_1- T_{M}(\zeta_1) \|_{N, \infty} S_{N,p}(1+ \delta) \le \left( \frac{\alpha}{8}\right)^{\frac{1}{p}}.
\end{align}
We set
\begin{align}
	L^{p}= \vert H \vert^{p}+ \vert u\vert^{p}+ \varphi^{p}+ \phi^{p}.
\end{align}
From $(\ref{eq: 4.7})$ and using $(\ref{eq: 4.9})$-$(\ref{eq: 4.10})$, we get
\begin{equation}\label{4.12}
	\begin{aligned}
		\| \nabla u \vert T_{K} u\vert^{\delta} \|_{p} \le C \| L \vert T_{K} u\vert^{\delta} \|_{p}.
	\end{aligned}
\end{equation}
We apply the H\"older inequality on the right-hand side of \eqref{4.12}. This leads to
\begin{equation}\label{4.14}
	\begin{aligned}
			\| \nabla u \vert T_{K} u\vert^{\delta} \|_{p} &\le C \| L \|_{r} \| T_{K}(u)\|_{\delta \frac{rp}{r-p}}^{\delta}\\
			& \le C  \| L \|_{r} \| \nabla \vert T_{K}(u) \vert^{\frac{r^{*}}{p^{*}}}\|_{p}^{ \delta \frac{p^{*}}{r^{*}} }\\
			& \le C \| L\|_{r} \| \nabla u \vert T_{K}(u) \vert^{\delta}\|_{p}^{\frac{\delta}{\delta+1}}.
	\end{aligned}
\end{equation}
Merging in \eqref{4.14}, we get
\begin{align}\label{4.15}
\| \nabla u \vert T_{K}(u) \vert^{\delta}\|_{p}^{\frac{p^{*}}{r^{*}}} \le C \| L \|_{r}.
\end{align}
Letting $K \rightarrow +\infty$ and taking into account the esitmate \eqref{4.15}, we have
\begin{equation}
	\begin{aligned}
		\| \nabla u \vert u \vert^{\delta}\|_{p}^{\frac{p^{*}}{r^{*}}} \le C \Big( \| H \|_{r} + \| \varphi \|_{r}+ \| \phi \|_{r} \Big),
	\end{aligned}
\end{equation}
i.e.
\begin{equation}
	\begin{aligned}
		\| \nabla \vert u \vert^{\frac{r^{*}}{p^{*}} }\|_{p} \le C \Big( \| H \|_{r} + \| \varphi \|_{r}+ \| \phi \|_{r} \Big)^{ \frac{r^{*}}{p^{*}}}.
	\end{aligned}
\end{equation}
Therefore, as long as $u \in L^{r}(\Omega)$ implies  $ \vert u \vert^{\frac{r^{*}}{p^{*}}} \in W_0^{1,p}(\Omega)$. This concludes our proof.
Clearly, the above argument works directly in the case $r< p^{*}$.  For the case $r>p^{*}$, we use a bootstrap procedure.	
\end{proof}


\begin{thebibliography}{9}
	\bibliographystyle{apa}
	
	\bibitem{Benilan} P. B\'enilan, L. Boccardo, T. Gallou\"et, R. Gariepy, M. Pierre, J. -L. Vazquez, \emph{An $L^1$- theory of existence and uniqueness of solutions of nonlinear elliptic equations}, Ann. Scuola Norm. Sup. Pisa 22 (1995) 241-273.
	
   \bibitem{Betta} M. F. Betta, A. Mercado, F. Murat, M. M. Porzio, \emph{Existence of renormalized solution to nonlinear elliptic equations with a lower-order term and right-hand side a measure}, J. Math. Pures Appl. 82 (2003) 90-124.
   
    \bibitem{Boca} L. Boccardo, \emph{Finite energy solutions of nonlinear Dirichlet problems with discontinuous coefficients}, 
   Boll. Un. Mat. It. 5 (2012) 357-368.
   
   
   
	\bibitem{Boccardo1} L. Boccardo, D. Giachetti, J. -I. Diaz, F. Murat, \emph{Existence and regularity of renormalized solutions to some elliptic problems involving derivations of nonlinear terms}, J. Differential Equations. 106 (1993) 215-237.
    
   	\bibitem{Boccardo } L. Boccardo, T. Gallou\"et, \emph{Nonlinear elliptic equations with right-hand side measure}, Comm. Partial Differential Equations 17 (1992) 641-955.
   	
   	\bibitem{Browder } F. E. Browder, \emph{Existence theorems for nonlinear partial differential equations}, Proceedings of Symposia in Pure Mathematics, Vol. 16, S. S. CHERN and S. SMALE Eds., A.M.S., Providence, 1970, 1-60.
   	
   	\bibitem{Carozza} M. Carozza, C. Sbordone,  \emph{The distance to $L^{\infty}$ in some function spaces and applications}, Differential Integral Equations 10 (1997) 599-607.
   	
   	\bibitem{Cirmi} G. R. Cirmi, S. D'Asero, S. Leonardi, M. M. Porzio, \emph{Local regularity results for solutions of linear elliptic equations with drift term}, Adv. Calc. Var. 15 1 (2022) 19-32.
   	
   	\bibitem{Farroni } F. Farroni, L. Greco, G. Moscariello, \emph{Noncoercive quasilinear elliptic operators with singular lower-order terms},  Calc. Var. 60 83 (2021).
   	
   	\bibitem{Giaquinta} M. Giaquinta, E. Giusti,  \emph{Quasi-minima}, Ann. Inst. H. Poincar\'re Anal. Non Linaire 1 2 (1984) 79-107.
   	
   	 	\bibitem{Gilbarg} D. Gilbarg, N.S. Trudinger, \emph{Elliptic partial differential equations of second order}, 249. Springer, New York, (2014).
   	
   	 \bibitem{Giannetti} F. Giannetti, L. Greco, G. Moscariello, \emph{Linear elliptic equations with lower-order terms}, Differential Integral Equations 26 5-6 (2013) 623-638.
   	 
   	  \bibitem{Guibe} O. Guib\'e, A. Mercaldo, \emph{Existence of renormalized solutions to nonlinear elliptic equations with two lower-order terms and measure data}, Trans. Amer. Math. Soc. 360 2 (2008) 643–669. 
   	 	 
   	\bibitem{GreMos} L. Greco, G. Moscariello, G. Zecca, \emph{Regularity for solutions to nonlinear elliptic equations}, Differential and Integral Equations 26 9-10 (2013) 1005-1113.
	
	\bibitem{LG} L. Greco, G. Moscariello, \emph{An embedding theorem in Lorentz--Zygmund spaces}, Potential Anal. 5 (1996) 581-590.
	\bibitem{LGreco} L. Greco, G. Moscariello, G. Zecca, \emph{An obstacle problem for noncoercive operators}, Abstr. Appl. Anal. (2015) 1-8.
	\bibitem{Greco} L. Greco, G. Moscariello, G. Zecca, \emph{Very weak solutions to elliptic equations with singular convection term}, J. Math. Anal. 457 2 (2018) 1376-1387.

	\bibitem{Lorentz} G. G. Lorentz, \emph{Some new function spaces}, Ann. Math. 51 (1950) 37-55.
	\bibitem{Moscariello} G. Moscariello, \emph{Existence and uniqueness for elliptic equations with lower-order terms}, Adv. Calc. Var. 4 4 (2011) 421-444.
   
   
   	\bibitem{Kim } H. Kim, T. P. Tsai, \emph{Existence, uniqueness, regularity results for elliptic equation with drift terms in critical weak spaces}, SIAM J. Math. Anal. 52 2 (2020) 1146-1191.
   
    \bibitem{ Lions} J. Leray, J. L. Lions, \emph{Quelques r\'esultats de Visik sur le probl\`emes elliptiques non lin\'eaires par les m\'ethodes de Minty-Browder}, Bull. Soc. Math. France 93 (1965) 97-107.
   
   	\bibitem{Maso } G. D. Maso, F. Murat, L. Orsina, A. Prignet, \emph{Renormalised solutions for elliptic equations with general measure data}, Ann. Scuola Norm. Sup. Pisa Cl. Sci. 28 (1999) 741-808.


	\bibitem{O'Neil} R. O'Neil, \emph{Convolution operators and $L(p,q)$ spaces}, Duke math. J. 30 (1963) 129-142.
	
	
	
	
	
	\bibitem{Porettea A} A. Porretta, \emph{Existence results for nonlinear parabolic equations via strong convergence of truncations}, Ann. Mat. Pura Appl. 177 4 (1999) 143-172.
	\bibitem{A. Porretta} A. Porretta, \emph{Weak solutions to Fokker--Planck equations and mean field games}, Arch. Ration. Mech. Anal. 216 1 (2015) 1–62.
	
	\bibitem{M.Pozrio} M. M. Pozrio, \emph{Existence, uniqueness and behavior of solutions for a class of nonlinear parabolic problems}, Nonlinear Anal. TMA 74 (2011) 5359-5382.
	\bibitem{Pozrio} M. M. Pozrio, \emph{On uniform and decay estimates for unbounded solutions of partial differential equations}, J. Diffential Equations 259 (2015) 6960-7011.
	
	\bibitem{Zadice} T. Radice, G. Zecca, \emph{Existence and uniqueness for elliptic equations with unbounded coefficients}, Ricerche Mat. 63 2 (2014) 355-367.
	
	\bibitem{Stampacchia} G. Stampacchia., \emph{Le probl\'eme de Dirichlet pour les \'equations elliptiques du second ordre \`a coefficients discontinuous}, Ann. Inst. Fourier (Grenoble) 15 (1965) 189-258.
	
	\bibitem{Zecca} G. Zecca, \emph{Existence and uniqueness for nonlinear elliptic equations with lower-order terms}, Nonlinear Anal. Theory Methods Appl. 75 (2012) 899-912.
	
	

	

\end{thebibliography}

\end{document}